\def\real{{\tt I\kern-.2em{R}}}
\def\Hyper#1{\hyper {\eskip #1}}%This is the change where the worthless Theorem 6.2 in section 6 is removed by removing section 6. The last revisions appear on both arxiv.org and vixra.org
\def\eskip{\hskip.25em\relax}
\def\nat{{\tt I\kern-.2em{N}}}
\def\hyper#1{\ ^*\kern-.2em{#1}}
\def\Hyper#1{\hyper {\eskip #1}}

\def\power#1{{{\cal P}(#1)}}
\def\qed{{\vrule height6pt width3pt depth2pt}\par\medskip}

\def\r#1{{\rm #1}}
\def\b#1{{\bf #1}}

\def\parm{\par\medskip}
\def\pars{\par\smallskip}
\def\m@th{\mathsurround=0pt}

\def\id{\par\hangindent2\parindent\textindent}
\def\textindent#1{\indent\llap{#1}}
\magnification=\magstep1
\tolerance 10000
\font\eightrm=cmr9
\baselineskip  14pt
\hoffset=.12in
\hsize 6.25 true in
\vsize 8.85 true in
\centerline{\bf General Logic-Systems that Determine Significant}
\centerline{\bf Collections of Consequence Operators}\par\medskip 
\centerline{Robert A. Herrmann}
\centerline{23 MAR 2006. Last revision 4 SEP 2013.}\par\bigskip
{\leftskip=0.5in \rightskip=0.5in \noindent {\eightrm {\it Abstract:} It is demonstrated how useful it is to utilize general logic-systems to investigate finite consequence operators (operations). Among many other examples relative to a lattice of finite consequence operators, a general logic-system characterization for the lattice-theoretic supremum of a nonempty collection of finite consequence operators is given.  Further, it is shown that for any denumerable language $L$ there is a rather simple collection of finite consequence operators and, for a propositional language, three simple modifications to the finitary rules of inference that demonstrate that the lattice of finite consequence operators is not meet-complete.  This also demonstrates that simple properties for such operators can be language specific. Using general logic-systems, it is further shown that the set of all finite consequence operators defined on $L$ has the power of the continuum and each finite consequence operator is generated by denumerably many general logic-systems.  Examples are given that define operators in terms of general logic-systems so that the physical entities produced require that the basic logic-system algorithm be applied. \par}}\pars Mathematics Subject Classifications (2000). 03B22, 03B65.\pars
 Keywords. Universal logic, general logic-systems, rules of inference, consequence operators. \par\bigskip
\noindent {\bf 1. Introduction.}\parm 
 In order to avoid an ambiguous definition for the ``finite consequence operator,'' it is assumed that a language $L$ is a nonempty set within informal set-theory (ZF). In the ordinary sense, a set $A \subset L$ is {\it finite} if and only if $A=\emptyset$ or there exists a bijection $f\colon A \to [1,n] = \{x\mid (1\leq x\leq n)\ {\rm and}\ (n \in \nat)\}$, where $\nat$ is the set of all natural numbers including zero. It is always assumed that $A$ is finite if and only if $A$ is Dedekind-finite. Finite always implies, in ZF, Dedekind-finite. There is a model $\eta$ for ZF that contains a set that is infinite and Dedekind-finite (Jech, 1971, pp. 116-118). On the other hand, for ZF, if $A$ is well-ordered or denumerable, then each $B\subset A$ is finite if and only if $B$ is Dedekind-finite. In all cases, if the Axiom of Choice is adjoined to the ZF axioms, finite is equivalent to Dedekind-finite. The definition of the general and finite consequence operator is well know but can be found in Herrmann (2006, 2004, 2001, 1987). \pars 

The subset map being consider has been termed as a (unary) ``operation.'' It has also termed either as a {\it consequence} or a {\it closure operator} by W\'ojcicki (1981). Due to its changed properties when embedded into a nonstandard structure, where for infinite $L$ the nonstandard extension of such a map is not a map on a power set to a power set but remains, at least, a closure operator, these two names were later combined to form the term {\it consequence operator} (Herrmann (1987)). In order to differentiate between two types, either the word {\it general} or {\it finite} ({\it or finitary}) is often adjoined to this term (Herrmann (2004)). Although finite consequence operators are closure operators with a finite character, they have additional properties, due to their set-theoretic definition, not shared, in general, by closure operators. Indeed, they have properties apparently dependent upon the construction of the language elements (Tarski, 1956, p. 71).   \pars

Since Tarski's introduction of consequence operator (Tarski, 1956, p. 60), although he mentions that it is not required for his investigations, a language $L$ upon which such operators are defined has been assumed to have, at the least, a certain amount of structure. For example, without further consideration, it has been assumed that $L$ can, at least, be considered as a semigroup or, often, a free algebra. Indeed, such structures have become  ``self-evident'' hypotheses. In order to emphasize that such special structures should not be assumed, the term ``non-organized'' is introduced (Herrmann (2006)). Although independent structural properties may exist, they are not considered in any manner as part of the hypotheses.\pars
 Formally, a {\it non-organized} $L$ is a language where only ``specifically stated'' properties $P_1,P_2, \ldots$ are assumed and where  either informal set theory or, if necessary, informal set theory with the Axiom of Choice is used to establish theorems informally. Hence, all other independent properties $L$ might possess are ignored. Indeed, the only property $L$ is assumed to possess is the method of ``word'' formation from a non-empty alphabet of symbols, images and other symbolized sensory information. When appropriate, the term ``non-specialized'' is only used as a means to stress this standard methodology.\parm

\noindent{\bf 2. General Logic-Systems.}\parm

In Herrmann (2006), the notion of a ``logic-system'' is discussed and an algorithm is described not in complete detail. The algorithm is presented here, in detail,  since it is  applied to most of the examples. In what follows, the algorithm, with associated objects, defines a {\it general logic-system} that when applied to a specific case yields {\it general logic-system deduction}.  The process is exactly the same as used in formal logic except for the use of the $RI(L)$ as defined below. Informally,  the pre-axioms is a nonempty $A\subset L$. (The term ``per-axioms'' is used so as not to confuse these objects with the notion of the ``consequence operator axioms'' $C(\emptyset)$.) The set of pre-axioms may contain any logical axiom and, in order not to include them with every set of hypotheses, $A$ can contain other objects $N \subset L$ that are consider as ``Theory Axioms'' such as natural laws as used for physical theories. There have been some rather nonspecific definitions for the rules of inference and how they are applied. It is shown in Herrmann (2006) that, for finite  consequence operators, more specific definitions are required. A {\it finitary rules of inference} is a fixed finite set ${RI(L)} =\{R_1,\ldots, R_p\}$  of $n$-ary relations $(0 < n \in \nat)$ on ${  L}.$ Note: it can happen that $RI(L) = \{\emptyset\}.$ (This corrects a misstatement made in Herrmann (2006, p. 202.) The pre-axioms are considered as a unary relation in $RI(L).$ An {\it infinite rules of inference} is a fixed infinite set ${RI(L)}$ of such $n$-ary relations on ${  L}.$  A {\it general rules of inference} is either a fixed finitary or infinite set of rules of inference. It is shown in Herrmann (2006), that there are finite consequence operators that require an infinite $RI(L)$, while others only require finite $RI(L).$  The term ``fixed'' means that no member of $RI(L)$ is altered by any set $X\subset L$ of hypotheses that are used as discussed below. All $RI(L)$, in this paper, are fixed.  For the algorithm, it is always assumed that an activity called {\it deduction} from a set of hypotheses $X \subset L$ can be represented by a finite (partial) sequence of numbered (in order) steps $\{b_1,\ldots,  b_{m}\}$  with the final step $  b_{m}$ a consequence (result) of the deduction. Also, $b_m$ is said to be ``deduced'' from $X.$ All of these steps are considered as represented by objects in the language $L.$ Each such deduction is composed either of the zero step,  indicating that there are no steps in the sequence,  or one or more steps with the last numbered step being some $  m >0$. In this inductive step-by-step construction,  a basic rule used to construct a deduction is the {\it insertion} rule. If the construction is at the step number $  m \geq 0,$ then the insertion rule,  {\bf I}, can be applied. This rule states: {\it Insertion of  any hypothesis (premise) from $  X \subset {  L},$ or insertion of a member from the set $  A,$ or the insertion of any member of any other unary relation can be made and this insertion is denoted by the next step number.} Having more than one unary relation is often very convenient in locating particular types of insertions. The pre-axioms are often partitioned into, at the least, two unary relations.  If the construction is at the step number $ m > 0,$ then ${RI(L)}$ allows for an additional insertion of a member from $  L$ as a step number $  m+1,$ in the following manner. For each $(j+1)$-ary $R_i,\ j \geq 1,$ if $  f \in   R_{i}$ and $  f(k) \in \{  b_1,\ldots,   b_{m}\},\   k=1,\ldots,  j,$ then $  f(j+1)$ can be inserted as a step number $  m+1.$ In terms of the notation $\vdash,$ where for $A \subset L$, $X\vdash A$ signifies that each $x\in A$ is obtained from some finite $F\subset X$ by means of a deduction, it follows from the above defined process that if $X\vdash b$, then there is either (1) a nonempty finite $F=\{b^1,\ldots,b^k\} \subset X$ such that $F \vdash b$ and each member of $F$ is utilized in $RI(L)$ to deduce $b$, or (2) $b$ is obtained by insertion of any member from any unary relation, or (3) $b$ is obtained using (2) by finitely many insertions and finitely many applications of the other $n$-ary ($n >1$) rules of inference. Hence, it follows that this algorithm yields the same ``deduction from hypotheses'' transitive property, as does formal logic, in that $X\vdash Y\subset L$ and $Y\vdash Z \subset L$ imply that $X\vdash Z.$\parm 

Note the possible existence of  special binary styled relations ${\bf J^\prime}$ that can be members of various $RI(L)$. These relations are identity styled relations in that the first and second coordinates are identical except that the second coordinate can carry one additional symbol that is fixed for the language used. In scientific theory building, these are used to indicate that a particular set of natural laws or processes does not alter a particular premise that describes a natural-system characteristic. The characteristic represented by this premise carries the special symbol and remains part of the final conclusion. Scientifically, this can be a significant fact. The addition of this one special symbol eliminates the need for the extended realism relation (Herrmann (2001)). Other deductions deemed as extraneous are removed by restricting the language. The deduction is constructed only from either the rule of insertion or the rules of inference via $AG$ (notation for the entire algorithm as described in this and the previous paragraph.) This concludes the definition of the logic-system. If $RI(L)$ is known to be either finitary or infinite, then the term ``general'' is often replaced by the corresponding term finite or infinite, respectively. 
\pars

For $ L,$ $  X \subset L$, general rules of inference $RI(L),$ and applications of $AG,$ the notation $RI(L) \Rightarrow C$ means that the map $C\colon \power {L} \to \power {L}$ (${\cal P}(L)=$ the power set of $L$) is defined by letting $C(X)= \{x\mid (X\vdash x)\ {\rm and}\ (x \in L)\}.$ The following result is established here not because its ``proof'' is complex, but, rather, due to its significance. Moreover, in Herrmann (2001), it is established in a slightly different manner and the result as stated there is not raised to the level of a numbered theorem.  Similar theorems relative to general consequence operators viewed as  closure operators have been established in different ways using a vague notion of deduction. What follows is a basic proof for the finite consequence operator using the required detailed definition for a general logic-system deduction. \parm
 {\bf Theorem 2.1} {\it Given non-specialized $L$, a general rules of inference ${RI(L)}$ and that the general logic-system algorithm $AG$ is applied. If $RI(L) \Rightarrow C$, then $C$ is a finite consequence operator.}\pars 

Proof. Let $C\colon{\cal P}(L) \to {\cal P}(L)$ be defined by application of the general logic-system algorithm ${AG}$ to each $X\subset L$ using the general rules of inference ${RI(L)}$. Let $x \in X.$ By insertion, $\{x\} \vdash x.$   Hence,  $  X \subset  C(X).$ If $  X \subset Y \subset L$ and $x\in C(X)$, 
then there is an $F \in {\cal F}(X)$ (= the set of all finite subsets of $X$) (= the set of all finite subsets of $X$)such that $F\vdash x$ and $F \subset Y.$ Hence, $x \in C(Y).$  Consequently,  $C(X) \subset C(Y).$  Let $y \in C(C(X)).$  From the definition of $C$, (1)  $X\vdash y$ if and only if $y \in C(X).$ By the transitive property for $\vdash,$ $C(X) \vdash C(C(X))$ implies  that $X\vdash C(C(X)),$ and (1) still holds. Hence, if $y \in C(C(X)),$ then $X \vdash y$ implies that $y \in C(X).$ Thus, $C(C(X))\subset C(X)$. Therefore, $C(X) = C(C(X))$ and $C$ is a general consequence operator. Let $x \in C(X).$ Then, as before, there is an $F {\cal F}(X)$ such that $F\vdash x$. Consequently, $C(X) \subset \bigcup\{C(F)\mid F\in {\cal F}(X)\} \subset C(X)$  and $C$ is a finite consequence operator. \qed

Let ${\cal C}_f(L)$ be the set of all finite consequence operators defined on $\power {L}.$ Each $C \in {\cal C}_f(L)$ defines a specific general rules of inference $RI^*(C)$ such that $RI^*(C) \Rightarrow C^* = C$ (Herrmann (2006)). However, in general, $RI(L)\not= RI^*(C)$.\pars

Let ${\cal C}(L)$ be the set of all general consequence operators defined on $\power {L}.$ Define on ${\cal C}(L)$ a partial order $\leq$ as follows: for $C_1,\ C_2 \in {\cal C}(L),$  $C_1 \leq C_2$ if and only if, for each $X \subset L,$ $C_1(X) \subset C_2(X).$ The structure $\langle {\cal C}(L),\leq  \rangle$ is a complete lattice. The meet, $\wedge$, is defined as follows: $C_1 \wedge C_2 = C_3,$ where for each $X \subset L,\ C_3(X) = C_1(X) \cap C_2(X).$ For each nonempty ${\cal H} \subset {\cal C}(L),\ \bigwedge {\cal H}$ means that, for each $X \subset L,$ $(\bigwedge {\cal H})(X) = \bigcap\{C(X)\mid C \in {\cal H}\}$ and, further, $\bigwedge {\cal H}=\inf\,{\cal H}.$  \pars

As is customary, in all of the following examples, explicit $n$-ary relations are represented in $n$-tuple form. Relative to the operator $\cup$, in the same manner as done in Herrmann (2006), if $\{a,b,c,d\} \subset L$, $\{\{(a,b),(c,d)\}\} \Rightarrow B,$ and $\{\{(a,c)\}\} \Rightarrow R,$ then defining $B \vee R$ as $(B \vee R)(X) = B(X) \cup R(X)= K(X)$ yields that $K \notin {\cal C}(L)$. Thus, ${\cal C}(L)$ is not closed under the $\vee$ operator as defined in this manner. Hence, if ``combined'' deduction is defined by this particular $\vee,$ then, in general, the combination does not follow the usual deductive procedures used through out mathematics and the physical sciences.\pars

Lemma 2.7 in Herrmann (2004) can be improved by simply assuming that ${\cal B} \subset {\cal P}(L),$ $L \in {\cal B}.$ The same proof as lemma 2.7 yields that the map defined by $C(X)= \bigcap \{Y\mid (X \subset Y)\ {\rm and}\  (Y \in {\cal  B})\} \in {\cal  C}(L).$ For a given $C \in {\cal  C}(L)$, $Y\subset L$ is a C-system ({\it closed system}) if and only if $Y = C(Y)$ (a closure operator fixed point). For each $C \in {\cal C}(L),$ let ${\cal S}(C)$ be the set of all C-systems. The equationally defined ${\cal  S}(C)= \{C(X)\mid X \subset L\}$ and $L \in {\cal S}(C).$ (If ${\cal B}$ is a {\it closure system} (i.e. closed under arbitrary intersection W\'ojcicki (1981) and $\cal B$ defines $C$, then ${\cal B} = {\cal S}(C).$) For nonempty ${\cal  H}\subset {\cal C}_f(L),$ let nonempty ${\cal S}' = \bigcap \{{\cal  S}(C)\mid C \in {\cal  H}\}.$ Using ${\cal  B} ={\cal  S}'$, if, for each $X \subset L,$ $(\bigvee_w{\cal  H})(X)= \bigcap \{Y\mid (Y\subset L)\ {\rm and}\  (X\subset Y)\ {\rm and}\  (Y \in {\cal  S}')\},$ then, for $\langle {\cal C}(L), \leq \rangle,\  \bigvee_w {\cal  H}=\sup\, {\cal  H}.$ The set of all consequence operators defined on ${\cal P}(L)$ forms a complete lattice $\langle {\cal  C}(L),\wedge,\vee_w,I,U\rangle$ with lower unit $I,$ the identity map, and upper unit $U,$ where for each $X \subset L,\ U(X) = L.$ If ${\cal  C}_f(L)$ is restricted to $\langle {\cal  C}(L),\wedge,\vee_w,I,U\rangle,$ then  $\langle {\cal  C}_f(L),\wedge,\vee_w,I,U\rangle$ is a sublattice. It is shown in Herrmann (2004), that $\langle {\cal  C}_f(L),\wedge,\vee_w,I,U\rangle$ is a join-complete sublattice.  (Note: Corollary 2.11 in the published version of Herrmann (2004) should read $\emptyset \not= {\cal A} \subset {\cal C}_f.$) Using finitary rules of inference, the fact that $\cup$ is not, in general, a satisfactory join operator for  $\langle {\cal S}(C), \subset \rangle$ is easily established. Consider non-specialized $L$ such that $\{a,b,c,d\} \subset L.$ Define $RI(L) = \{\{(a,c)\},\{(a,b,c,d)\}\} \Rightarrow B.$ Then $B(\{b\}) \cup B(\{a\}) = \{a,b,c\}.$ But, $\{a,b,c\}$ is not a C-system for $B$ since $B(\{a,b,c\}) = \{a,b,c,d\}.$  Defining for each $C \in {\cal C}(L)$ and each $X,Y \in {\cal S}(C),\ X\uplus Y = C(X\cup Y),$ then the structure $\langle {\cal S}(C), \subset \rangle$ is a complete lattice with the join $\uplus$ and meet $X \wedge Y = X \cap Y.$ 
\pars

For each non-specialized language $L$ and non-empty ${\cal H} \subset {\cal C}_f(L),$ a natural investigation would be to determine whether there is a significant relation between $\bigvee_w{\cal H}$ and any collection of general logic-systems that generates each member of $\cal H.$ For each $C \in {\cal H},$ let $RI_C(L)$ be any general rules of inference such that $RI_C(L) \Rightarrow C$. \parm

{\bf Theorem 2.2.} {\it If $L$ is non-specialized, then for the structure $\langle {\cal C}_f(L),\wedge,\vee_w,I,U\rangle$ and each  nonempty ${\cal  H} \subset {\cal  C}_f(L),$ it follows  that 
$\bigcup\{RI_C(L)\mid C \in {\cal H}\} \Rightarrow \bigvee_w{\cal H}.$}\pars 
Proof. For ${\cal H},$ let $\bigcup\{RI_x(L)\mid x \in {\cal H}\} \Rightarrow {\cal U},$ $X \subset L,$  and $C \in {\cal H}.$ Since $C \leq {\cal U},$ then ${\cal U}(X) \subset C({\cal U}(X))\subset {\cal U}({\cal U}(X))= {\cal U}(X)$ implies that ${\cal U}(X) = C({\cal U}(X)).$ Thus, for each $C \in {\cal H}$, ${\cal U}(X)$ is a C-system and, hence, ${\cal U}(X) \in {\cal S}'= \bigcap \{{\cal S}(C)\mid C \in {\cal H}\}.$\pars

Suppose that $X \subset Y \in {\cal S}'.$ Then, for each $C\in {\cal H},\ X \subset Y =C(Y)$ implies that, for each $C\in {\cal H}, \ X \subset {\cal U}(X) \subset {\cal U}(C(Y))$. Consider $b \in {\cal U}(C(Y)).$ Take any finite $F \subset Y=C(Y)$ such that $F$ is used to obtain $b$ by application of $AG$ as the next step in a deduction using $\bigcup\{RI_x(L)\mid x \in {\cal H}\}$. Then $F$ is used along with finitely many ($\geq 0$) $RI_{C_i}(L)\Rightarrow C_i \in {\cal H}$ to obtain $\{b_1,\ldots,b_m\}.$ Since for each $i \in [1,k],\ b_i \in C'(Y) = Y,$ for some $C' \in {\cal H},$ then $\{b_1,\ldots,b_m\}\subset Y.$  If $b\notin \{b_1,\ldots,b_n\},$ then there are finitely many ($\geq 0$) $RI_{C_j}(L) \Rightarrow C_j \in {\cal H}$ and from $F$ and $\{b_1,\ldots,b_n\}$ the set $\{c_1,\ldots,c_k\}$ is deduced. But again $\{c_1,\ldots,c_k\} \subset Y.$ This process will continue no more than finitely many times until $b$ is obtain as a member of a finite set of deductions from members of $\bigcup\{RI_x(L)\mid x \in {\cal H}\}$ and $b \in Y.$ Hence, ${\cal U}(C(Y))\subset Y.$
  But, $C(Y) =Y$ implies that $Y \subset {\cal U}(C(Y)).$ Hence, $Y = {\cal U}(C(Y)) = {\cal U}(Y)$ and, since ${\cal U}(X) \subset {\cal U}(Y)$, then ${\cal U}(X) \subset Y=C(Y)$ for each $C\in {\cal H}.$ Therefore, ${\cal U}(X) \subset Y \in {\cal S}'.$ Hence, ${\cal U}(X) = (\bigvee_w{\cal H})(X)$. \qed

After showing that ${\cal C}_f(L)$ is closed under finite $\wedge$, then Theorem 2.2 yields a general logic-system proof that $\langle {\cal C}_f(L),\wedge,\vee_w,I,U\rangle$ is a join-complete lattice. It is rather obvious that, in general, if $RI_C(L) \Rightarrow C$ and $RI_D(L) \Rightarrow D,$ then $RI_C(L) \cap RI_D(L) \not\Rightarrow C \wedge D.$ For example, let $\{a,b,c,d\} \subset L$ and $RI_C(L) = \{\{(a,b)\}\},\ RI_D(L) = \{\{(a,b),(b,c)\}\}.$ Then $C(\{a\}) = \{a,b\},\ D(\{a\}) = \{a,b,c\}$ implies that $(C\wedge D)(\{a\}) = \{a,b\}.$ But, $RI_C(L) \cap RI_D(L) = \emptyset \Rightarrow I$ and $I(\{a\}) = \{a\}.$ Even if we took the intersection, $\cap_1,$ of the individual relations from each general rules of inference, then, for $RI_E(L) = \{\{(a,b),(b,c)\}\}$ and $RI_F(L) = \{\{(a,b),(b,d),(d,c)\}\},$ it would follow that $RI_E(L) \cap_1 RI_F(L) \not\Rightarrow E \wedge F.$ However, it is obvious that, for each nonempty ${\cal H} \subset {\cal C}_f(L),$ if $\bigcap\{RI_x(L)\mid x \in {\cal H}\} \Rightarrow G \in {\cal H},$ then $G = \bigwedge\, {\cal H}.$\pars

There is a constraint that can be placed on deduction from hypotheses using algorithm $AG$. With one exception, there is a $RI(L)$ that if the restricted $RI(L) \Rightarrow D$, then $D$ is not a general consequence operator. 
 \parm

{\bf Example 2.2.} ({\it Limiting the number of steps in an $RI(L)$-deduction need not yield a consequence operator.}) Suppose that $AG$ has the added restriction that no deduction from hypotheses be longer then $n$ steps, where $n>1.$ For each $L,$ such that $\vert L\vert \geq n+1,$ let $a \not= b,$ for $i \in [1,n-1],\ x_i \notin \{a,b\},\ \{x_i,a,b\} \subset L,$ and if $i,j \in [1,n-1], \ i \not= j,$ then $x_i \not= x_j.$ Consider $RI(L)= \{\{(x_1,\ldots,x_{n-1},a)\},\{(a,b)\}\}.$ Let $\vdash_{\leq n}$ indicate that each deduction from premises, using $RI(L),$ most have $n$ or fewer steps. Then, using this restriction, for $X \subset L,$ let $D(X) =\{x\mid (X\vdash_{\leq n} x)\ {\rm and}\  (x \in L)\}.$ Consider $X = \{x_1,\ldots,x_{n-1}\}$. Then $D(X)=X \cup \{a\}.$ But $D(D(X))= D(X \cup \{a\}) = X \cup \{a,b\}.$ This  follows since the definition requires that you calculate in no more than $n$ steps {\it all} of the consequences of $\{x_1,\ldots,x_{n-1}, a\}$ using {\it any} finite subset of $\{x_1,\ldots,x_{n-1}, a\}.$  Thus, $D^2 \not= D$ and $D \notin {\cal C}(L).$ Let $PR$ be a standard predicate language (Mendelson, 1987, pp. 55-56), where $PR$ has more than one predicate with one or more arguments and with the set of variables $\cal V$.  Let $R^1$ be the set of all axioms, $R^2 = \{(A,(\forall xA))\mid (x \in {\cal V})\ {\rm and}\ (A \in PR)\}$ and $R^3 = \{(A\to B),A,B)\mid A, B \in PR\}.$ If you restrict predicate deduction to 3 steps or less, then restricted $RI(PR) \Rightarrow C_P$ and $C_P$ is not a general consequence operator. \parm

\noindent {\bf 3. Special Consequence Operators.}\parm

Throughout this section, unless other specific properties are stated, the language $L$ is non-specialized. In Herrmann (1987), two significant collections of consequence operators are defined. Let $X\cup Y \subset L.$ (1) Define the map $C(X,Y)\colon \power {L}\to  \power {L}$ as follows: for $A \in \power {L}$ and $A \cap Y \not= \emptyset,\ C(X,Y)(A)= A \cup X$. If $A \cap Y = \emptyset,\ C(X,Y)(A)=A.$ (2) Define the map $C'(X,Y)\colon \power {L}\to  \power {L}$ as follows: for $A \in \power {L}$ and $Y \subset A,\ C'(X,Y)(A)= A \cup X$. If $Y\not\subset A,\ C'(X,Y)(A)=A.$ It is shown in Herrmann (1987) via long set-theoretic arguments that each $C(X,Y) \in {\cal C}_f(L),$ and $C'(X,Y) \in {\cal C}(L).$ If $Y \in {\cal F}(L),$ then $C'(X,Y) \in {\cal C}_f(L).$  Now suppose that $Y$ is infinite and $Y \subset A.$ Then for each $F \in {\cal F}(L)$, since $Y \not\subset F$, then $C'(X,Y)(F) = F.$ Hence, $\bigcup \{C'(X,Y)(F)\mid F\in {\cal F}(A)\} = A$. But if $X \not\subset A,$ then $C'(X,Y)(A) = A \cup X \not=
 \bigcup \{C'(X,Y)(F)\mid F\in {\cal F}(A)\}.$ Therefore, if infinite $Y\subset A\subset L,$ and $X\not\subset A,$ then $C'(X,Y) \in {\cal C}(L) -{\cal C}_f(L).$ Thus, in general, for infinite $L,\ C'(X,Y)$ need not be finite. \pars

 In some cases, the use of logic-systems can lead to rather short proofs for consequence operator properties, where other methods require substantial effort.\parm 

{\bf Example 3.1.} ({\it An obvious sufficient condition for $\bigwedge {\cal H} \in {\cal C}_f(L),$ when nonempty ${\cal H} \subset {\cal C}_f(L)$})  For non-specialized $L$, let  nonempty ${\cal H} \subset {\cal C}_f(L).$ If $\bigcap\{RI_x(L)\mid x \in {\cal H}\} \Rightarrow G \in {\cal H},$ then $G = \bigwedge\, {\cal H}.$ \qed\parm

{\bf Example 3.2.} ({\it Establishing that some significant general consequence operators are finite.}) We use logic-systems to show that $C(X,Y) \in {\cal C}_f(L)$ and, if $Y \in {\cal F}(L),\   X \subset L,$ then $C'(X,Y)$ is finite. For $C(X,Y)$ if $Y$ or $X =\emptyset,$ let $RI(L) =\emptyset\Rightarrow I.$ If $Y$ and $X \not=\emptyset,$ let $RI= \{R^2\},$ where $R^2 = \{(y,x)\mid (y\in Y)\ {\rm and}\  (x \in X)\}.$ Then it follows easily that $RI(L) \Rightarrow C(X,Y).$ Thus, $C(X,Y)$ is finite. If $X = \emptyset,$ then $C'(Y,X) = I$  and $RI'(L) = \emptyset \Rightarrow I.$ Now let $Y\in {\cal F}(L).$ If $Y =\emptyset$ and $X\not=\emptyset,$ then let $RI'(L) = \{R^1\},$ where $R^1= X.$ If $X$ and $Y \not= \emptyset,$ then there is an bijection $f\colon [1,n] \to Y$. In this case, let $RI'(L) = \{\{(f(1),\ldots,f(n),x)\mid x \in X\}\}.$ Then $RI'(L) \Rightarrow C'(X,Y).$ Hence, if $Y \in{\cal F}(L),$ then $C'(X,Y) \in {\cal C}_f(L).$   \qed\parm 

Relative to a standard propositional language $PD$, after some extensive analysis and using the \L o\'s and Suszko matrix theorem,  W\'ojcicki (1973) defines a collection of $k$-valued matrix generated finite consequence operators $\{C^*_k\mid k = 2,3,4,\ldots\}$ such that the greatest lower bound for this set in the lattice $\langle {\cal C}(PD), \leq \rangle$ is not a finite consequence operator. Are there simpler examples that lead to the same conclusion? \parm
 {\bf Example 3.3.} ({\it Showing that, in general, $\langle {\cal C}_f(L), \wedge, \vee_w, I,U \rangle$ is not a meet-complete lattice.}) Let $L$ be any denumerable language. Hence, there is a bijection $f\colon \nat \to L.$ Define $B_n= f[[1,n]]$ for each $n \in \nat^{>0},$ where $\nat^{>0} = \{ n\mid (n\in \nat)\ {\rm and}\  (n \geq 1)\}.$  Then for each $n \in \nat^{>0},\ f(0) \not\in B_n$. Let $X = \{f(0)\}$ and $C_n = C'(X, B_n).$ 
We have that $\inf \{C'(X, B_n) \mid (n\geq 1)\ {\rm and}\  (n \in \nat)\} = C'(X, f[\nat]-\{f(0)\})\leq C'(X,B_n)$ for each $B_n.$ But, since $f[\nat]-\{f(0)\}$ is an infinite set and, for  $A = f[\nat]-\{f(0)\}, X\not\subset A,$ then $C'(X, f[\nat]-\{f(0)\})$ is not a finite consequence operator. The fact that this consequence operator is not finite also holds for non-denumerable infinite $L,$ where $L$ either has additional structure, or an additional set-theoretical axiom such as the Axiom of Choice is utilized.  \qed\parm 

Of course, $C'(X,Y)$ is not the usual type of consequence operator one would associate with a propositional language. Are there simple finite consequence operators associated with standard formal propositional deduction that are not meet-complete? \pars

Using finite logic-systems, the following examples show how various weakenings for deduction relative to, at least, a propositional language $PD,$ generate collections of consequence operators that also establish that $\langle {\cal  C}_f(PD),\wedge,\vee_w,I,U\rangle$ is not a meet-complete lattice.\pars

The propositional language $PD$ defined by denumerably many (distinct) propositional variables $P=\{P_n\mid n \in {\nat}\},$ and is constructed in the usual manner from the unary $\neg$ and binary $\to$ operations. For the standard propositional calculus and deduction, one can use the following sets of axioms, with parenthesis suppression applied.
$R_1 = \{X \to (Y \to X)\mid (X \in PD)\ {\rm and}\ (Y \in PD)\},\ R_2 = \{(X \to (Y\to Z))\to ((X\to Y)\to (X\to Z)) \mid (X \in PD)\ {\rm and}\ (Y \in PD) \ {\rm and }\ (Z \in PD) \},\ R_3 = \{(\neg X\to \neg Y) \to (Y \to X) \mid (X \in PD)\ {\rm and}\ (Y \in PD)\}.$ The one rule of inference $MP = R^3(PD)=\{(X\to Y,X,Y)\mid (X \in PD)\ {\rm and}\ (Y \in PD)\}.$ Let $R^1(PD) = R_1 \cup R_2 \cup R_3.$ Standard proposition deduction $PD$ uses the rules of inference $RI(PD) = \{R^1(PD), R^3(PD)\}\Rightarrow C_{PD}.$ Let ${\cal T}$ be the set of all $PD$ tautologies under the standard valuation. Then by the soundness and completeness theorems ${\cal T} = C_{PD}(\emptyset).$ In all of the following examples, $R_1,\ R_2,\ R_3,\ R^1(PD),\ R^3(PD)$ are as defined in this paragraph and $RI(PD)$ is modified in various ways\parm

{\bf Example 3.3.1.} ({\it Propositional deduction with a restricted Modus Ponens rule yields $\{C_n\} \subset {\cal C}_f(L)$ such that $\bigwedge \{C_n\} \notin {\cal C}_f(L).$}) Consider $PD.$ Let ${\cal J} = \{((  P_i\to   P_0),P_i, P_0)\mid i \in \nat^{>0}\}.$ Let ${H} = R^3(PD) - {\cal J}.$ For each $n \in \nat^{>0},$ let $R^3_n = { H} \cup \{((  P_n\to   P_0), P_n, P_0)\}.$ Thus, the Modus Ponens rule of inference is restricted for each $n \in \nat^{>0}.$ Let $RI_n(PD)= \{R^1(PD), R^3_n\} \Rightarrow C_n.$ Now let $X = \{(  P_n\to   P_0),P_n \mid n\in \nat^{>0}\}.$ Then, for all $n\in \nat^{>0},\ P_0 \in C_n(X).$ Hence, $P_0 \in (\bigwedge \{C_n\})(X).$ Consider for any $n \in \nat^{>0},\ F \in {\cal F}(X)$ such that $P_0 \in C_n(F).$ Since $P_0 \notin {\cal T}$, then $P_0 \notin C_n(\emptyset)$ implies that $F \not=\emptyset.$ Further, for some $k \in \nat^{>0},\ \{(  P_k \to   P_0), P_k\} \subset F.$ For, assume not. First, consider, for $n \in \nat^{>0},$ $\{(  P_j \to   P_0), P_k \} \subset F,\ \{k,j\}\subset \nat^{>0},\ k \not= j$ and assume that $(  P_j \to   P_0), P_k \vdash_n P_0.$ This implies that $\vdash_n (  P_j \to   P_0) \to (P_k \to P_0),$ where the part of the Deduction Theorem being used here does not require any of the objects removed from the original $R^3(PD)$.  But, $\vdash_n $ implies $\models_{PD}$, using the standard valuation which is not dependent upon our restriction. Hence. $\models_{PD} (  P_j \to   P_0) \to (P_k \to P_0).$ However, $\not\models_{PD} (  P_j \to   P_0) \to (P_k \to P_0).$ The same would result, for $k \in \nat^{>0},$ if only the wwfs $P_k,$ or only wwfs $(  P_k \to   P_0)$ are members of $F$. Hence, there exists a unique $M = \max \{i \mid ((  P_i \to   P_0) \in F){\rm\ and}\ (P_i \in F)\ {\rm and}\  (i\in \nat^{>0})\}.$ But, then $P_0 \notin C_{M+1}(F).$ Consequently, 
this implies that $P_0 \notin (\bigwedge \{C_n\})(F).$ Thus, $\bigcup \{(\bigwedge \{C_n\})(F) \mid F\in {\cal F}(X)\}\not= (\bigwedge \{C_n\})(X)$ yields that $\bigwedge \{C_n\} \in {\cal C}(PD) -{\cal C}_f(PD).$ 
 \qed \parm 

For each $R\subset R^1(PD),$ always consider the standard elementary valuations for propositional wwfs. Also, if $R \subset R^1(PD),\ X \subset PD,$ and one considers the rules of inference $RI_R(PD) = \{R, R^3(PD)\}\Rightarrow C_R,$  then $X \vdash_R A$ implies that $X\vdash_{PD} A.$ Hence, if  $X \vdash_R A$, then, for each $x \in A,$ there is some $F\in {\cal F}(X)$ such that $F \models_{PD} x.$  Although, ${\cal T} = C_{PD}(\emptyset),$ in general, ${\cal T}\not= C_R(\emptyset).$  However, we do have that ${\cal T} \supset C_R(\emptyset).$    \parm

{\bf Example 3.3.2.} ({\it $PD$ axioms with a missing atom $P_0$ yields $\{C'_m\} \subset {\cal C}_f(PD)$ such that $\bigwedge \{C'_m\} \notin {\cal C}_f(PD).$})  Consider $PD.$ Let $L'$ be the propositional language defined by the set of propositional variables $\{P_i\mid i \in \nat\}-\{P_0\}.$ For each $m \in \nat^{>0},$ let $J_m =(\neg P_0 \to \neg P_m) \to (P_m \to P_0),$ and let $R_1',\ R_2',\ R_3'$ be defined for the language $L'$, in the same manner as $R_1,\ R_2,\ R_3$ are defined for $L,$ and let $R^3(PD)$ be defined for $PD.$ Let $R^1 = R_1'\cup R_2'\cup R_3',$ and, for each $m \in \nat^{>0},$ $R^1_m = \{R^1\cup \{J_m\}\}.$ For each $m \in \nat^{>0},$ the rules of inference is the set $RI_m'(PD) = \{R^1_m, R^3(PD)\} \Rightarrow C'_m$ and, for this rules of inference, the $P_0$ only appears in $J_m \cup R^3(PD).$ For any deduction, the Modus Ponens (MP) rule is applied to previous steps. Thus, no deduction, from empty hypotheses,using $R^1$ can either lead to any wwf that includes $P_0$ or utilize any wwf that contains $P_0$. The only member of the $R^1_m$ that is not a premise and can be used for a deduction that contains $P_0$ is $J_m$. Let $X = \{(\neg P_0 \to \neg P_n), P_n\mid n \in \nat^{>0}\}.$ Obviously, for each $m \in {\nat}^{>0},\ P_0 \in C'_m(X)$ and, since $J_m \in {\cal T}$ and $P_0 \notin {\cal T},$ then $P_0 \notin C'_m(\emptyset).$  Consider for each $m \in \nat^{>0},$ nonempty $A \in \{J_n, (\neg P_0 \to \neg P_n), P_n, P_0\mid (m \not= n\in \nat^{>0})\}$. Then $\not\vdash_m A.$ For example, let $A = J_n\ n\not=m.$ This would imply that $\vdash_m J_n.$ But, since $J_m \not= J_n$ and there is no member of $R^1$ to which $MP$ applies, such a deduction is not possible. The same holds for $(\neg P_0 \to \neg P_n), P_n, P_0.$ Further, for $A$ and for $j \not=m$ or $k \not=m$, $(\neg P_0 \to \neg P_j), \ \neg P_k\not\vdash_m P_0$ for the same reasons.  Consider for each $m \in \nat^{>0}$, any nonempty $F \in {\cal F}(X)$ such that $P_0 \in C'_m(F).$ Then, from the above discussion,  $(\neg P_0 \to \neg P_m), P_m \in F.$ Let $a = \max \{i \mid ((\neg P_0 \to P_i)\in F)\ {\rm and}\  (i\in \nat^{>0})\},\  b = \max \{i\mid (P_i \in F)\ {\rm and}\  (i\in \nat^{>0})\}.$ Let $M = \max\{a,b\}$. Then, again from the above discussion, $P_0 \notin C'_{M+1}(F).$ Hence, $P_0 \notin \bigcup \{(\bigwedge  \{C'_m\})(F)\mid F\in {\cal F}(X)\}\not= (\bigwedge \{C'_m\})(X)$ and $\bigwedge \{C'_m\} \in {\cal C}(PD) - {\cal C}_f(PD).$    \qed 
     
{\bf Example 3.3.3.} ({\it Extended positive propositional deduction ($PD$ axiom restrictions) yields $\{C_n\} \subset {\cal C}_f(L)$ such that $\bigwedge \{C_n\} \notin {\cal C}_f(L).$})  
Consider $PD$. As defined above ${\cal T}$ is the set of all $A\in PD$ such that $A$ is a tautology. The h-rule is defined as follows: for each $A \in L,$ let $h(A)$ denote the wwf that results from erasing each $\neg $ that appears in $A$. Now let $R_3' = \{X \mid (X \in R_3)\ {\rm and}\  (h(X) \in {\cal T})\}.$ Then $\emptyset \not= R_3'\not= R_3$ since if $h(A)\in {\cal T},$ then $h((\neg A \to \neg B) \to (B \to A)) = (h(A) \to h(B)) \to (h(B)\to h(A)) \in {\cal T}$ and $(\neg P_0 \to \neg P_n)\to (P_n \to P_0) \notin R'_3,\ n \not= 0.$ Let $R^1 = R_1 \cup R_2\cup R_3'$ and $RI_h(PD) = \{ R^1, R^3(PD)\}\Rightarrow C_h.$ For each $n \in \nat^{>0},$ let $J_n = (\neg P_0 \to \neg P_n) \to (P_n \to P_0)$ and the rules of inference be $RI_n(PD) = \{R^1\cup \{J_n\},R^3(PD)\} \Rightarrow C_n.$  Each member of $R^1$ is a tautology. Further, if $A \in R^1,$ $h(A) \in {\cal T}$ and if $A,\ A\to B \in R^1,$ then $h(A\to B) = h(A)\to h(B)$ implies that $h(B) \in {\cal T}.$ Thus, for each $A \in R^1$, the $h$ operator coupled with  any MP application using members of $R^1$ yields a tautology. This operator acts as a concrete model for deduction from empty hypotheses using  members of $R^1.$ But for certain members of $R_3,$ the h-rule does not generate a tautology and these members of $R_3$ are, therefore, not members of $C_h(\emptyset).$ That is, for $R_1 \cup R_2\cup R_3'$ they are not $RI_h(PD)$ theorems. Each $J_n$ is a wwf that cannot be established by $RI_h(PD)$ deduction (i.e. $J_n \notin C_h(\emptyset)$). Consider for any $n \in \nat^{>0}$,  $A\vdash_n B.$ This can always be written as $J_n, A \vdash_n B.$ 
Suppose that for each $m,n,k \in \nat^{>0},\ k \not=n,$ that $X_m = (\neg P_0 \to \neg P_m)$ and $X_m, P_k \vdash_n P_0.$ Since the derivation of the Deduction Theorem does not utilize $R_3$, then this implies that $\vdash_n J_n \to (X_m \to (P_k \to P_0)).$ This can be considered as a deduction that does not use $J_n$ as a premise. Hence, this implies that $\vdash_h J_n \to (X_m \to (P_k \to P_0)).$ However, this contradicts the h-rule. Also notice that $J_m = (X_m \to (P_m \to P_0)).$ 
Hence, for each $m,n,k \in \nat^{>0}, k\not= n$; $X_m, P_k \not\vdash_n P_0,$ implies that for any nonempty $A \subset \{X_m, P_k\mid m,k \in \nat^{>0})\ {\rm and}\ (k \not = n)\},$ that $P_0 \notin C_n(A).$ However, for each $n \in \nat^{>0},\ P_0 \in C_n(\{X_n,P_n\}).$  This also shows that for each $m,n \in \nat^{>0}, n \not =m,$ that $C_n(\{X_m,P_m\}) \not= C_m(\{X_m,P_m\}),$ and that $C_n \not= C_m.$ Obviously, since $P_0 \notin {\cal T}$ implies that, for each $n \in \nat^{>0},\ \not\vdash_n P_0$, then, for each $n \in \nat^{>0},\ P_0 \notin C_n(\emptyset).$ Now let $Y=\{(\neg P_0 \to \neg P_i),P_i\mid i \in \nat^{>0}\}.$ Then, for each $n \in \nat^{>0},\ P_0 \in C_n(Y).$ Thus $P_0 \in (\bigwedge \{C_n\mid n \in \nat^{>0}\})(Y).$ Consider for each $j \in \nat^{>0}$, any $F \in {\cal F}(Y)$ such that $P_0 \in C_j (F).$ Then $F \not=\emptyset.$ If $\{i\mid ((\neg P_0 \to \neg P_i)\in F)\ {\rm and}\  (i \in \nat^{>0})\}\not=\emptyset$, let $a=\max\{i\mid ((\neg P_0 \to \neg P_i)\in F){\ \rm and\ }(i\in \nat^{>0})\}$. If $\{i\mid (P_i\in F)\ {\rm and}\  (i \in \nat^{>0})\}\not=\emptyset,$ let $b = \max\{i\mid (P_i\in F)\ {\rm and}\  (i \in \nat^{>0})\}$. The set $\{a,b\}\not=\emptyset.$ Let $M=\max \{a,b\}.$ It has been shown that $P_0 \notin C_{M+1}(F).$  Hence, from this, it follows that $P_0 \notin \bigcup \{(\bigwedge  \{C_n\})(F)\mid F\in {\cal F}(Y)\}\not= (\bigwedge \{C_n\})(Y)$ and $\bigwedge \{C_n\} \in {\cal C}(PD) - {\cal C}_f(PD).$\qed

For the two collections $\{C_n\},\ \{C_m\} \subset {\cal C}_f(L)$ defined in the last two examples, notice that $\bigcap RI_m'(PD) = \bigcap RI_n(PD) = \{R^3(PD)\} \Rightarrow G \in {\cal C}_f(L),\ G(\emptyset) = \emptyset,\ G < \bigwedge \{C_n\}.$ The rule of inference $\{R^3I(PD)\}$ yields axiomless propositional deduction.\parm

{\bf Example 3.4.} ({\it For denumerable $L$, the set ${\cal C}_f(L)$ has the power of the continuum.})
For any set $X$, let $\vert X \vert$ denote its cardinality (power). For the real numbers $\real$, $\vert \real \vert$ is often denoted by $\aleph$ or $c$. For a denumerable language $L$, let $a \in L$ and consider $L - \{a\}.$ Let ${\cal I}$ be the set of all infinite subsets of $L - \{a\}.$ Then $\vert {\cal I} \vert = \aleph.$ For any $X \in {\cal I}$, let $R_X = \{(a,x)\mid x \in X\}$ and $RI_X(L)=\{R_X\}\Rightarrow C_X.$ Then $C_X(\{a\}) = \{a\}\cup X.$ Let $A,B \in {\cal I}, \ A \not= B.$ Then $C_A(\{a\}) = \{a\} \cup A \not= \{a\} \cup B = C_B(\{a\}).$ Thus $\vert \{C_X\mid X \in {\cal I}\} \vert = \aleph.$ Hence $\vert {\cal C}_f(L)\vert \geq \aleph.$ \pars
On the other hand, each $C \in {\cal C}_f(L)$ corresponds to a general logic-system $RI^*(C)$ such that $RI^*(C) \Rightarrow C$ (Herrmann (2006)). From the definition of a general rules of inference, $RI^*(C)$ corresponds to a finite or denumerable subset of $\bigcup (\{L^n\mid n \in \nat^{>0}\}$. But, $\power {\bigcup (\{L^n\mid n \in \nat^{>0}\}} = \aleph.$   
Hence, $\vert {\cal C}_f(L)\vert \leq \aleph.$ Consequently, $\vert {\cal C}_f(L)\vert =\aleph.$ (Depending upon the definition of ``infinite,'' this result may require the Axiom of Choice.) \qed

{\bf Example 3.5.} ({\it For denumerable $L,$ there exists denumerably many general logic-systems that generate a specific $C \in {\cal C}_f(L).$}) Let $C\in {\cal C}_f(L).$ Let $RI^*(C)$ be the general logic-system defined in Herrmann (2006), where $RI^*(C) \Rightarrow C.$ Notice that when the $RI^*(C)$-deduction algorithm is used, it can be considered as applied to $\bigcup RL^*(C)$. For $\emptyset \not= X \in {\cal F}(L)$, where $\vert X\vert = n \in \nat$ and $n\geq 1,$ consider any finite sequence $\{x_1,\ldots, x_n\} = X.$
Define $R_X = \{(x_1,\ldots,x_n,x)\mid x \in X\}.$ Let general logic-system $RI_1(L)= \{R_X\mid X \in {\cal F}(L)\}.$ Then $RI_1(L) \Rightarrow C_1 \in {\cal C}_f(L).$ Let $Y \in\power {L}.$ If $Y = \emptyset,$ then $C_1(\emptyset) = \emptyset.$ For nonempty $Y \in \power {L},$ let $y \in C_1(Y),$ then $y$ is deduced via the general logic-system algorithm. Hence, there exists a nonempty finite $A = \{y_1,\ldots,y_n\}= Y\subset L$ such  that $(y_1,\ldots,y_n,y) \in RI_1(L)$ and $y \in Y.$ Hence, $C_1(Y) \subset Y$ implies that $C_1(Y) = Y.$ Thus, $C_1$ is the identity finite consequence operator. \pars

Let $RI^+(L) = RI_1(L) \cup RI^*(C)$ and note that $RI^+(L) \Rightarrow C.$ For each $n \in \nat^{>0},$ there exists $r_n \in \bigcup RI^+(L),$ such that $r_n = (x_1,\ldots,x_n,x), \ i = 1,\ldots, n$ and $x \in C(\{x_1,\ldots,x_n\}).$ Thus, there exists a unique nonempty $R_n^+ \subset \bigcup RI^+(L)$ such that $r_n \in R_n^+$ if and only if $p_i(r_n) = x_i\in L, \ 1,\ldots n.$ The general logic-system $RI^{**}(L) = \{R^1\} \cup \{R^+_k\mid k\in \nat^{>0}\} \Rightarrow C,$ where $R^1 = C(\emptyset).$ (Notice that if $A \subset R^1,$ then $C(A) = R^1.$) For each $n \in \nat, \ n \geq 2,$ let $(y_1,\ldots,y_n)$ be a distinct permutation $p$ of the coordinates $x_i,\ i = 1,\ldots, n,$ for a specific $r_n = (x_1,\ldots, x_n,x) \in R^+_n.$ Let $r_n^p = (y_1,\ldots,y_n,x)$ and $R^+_{n,p} = (R^+_n - \{r_n\}) \cup \{r_n^p\}.$  This yields $RI^p_n(L) = (RI^{**}(L) - \{R^+_n\}) \cup \{R^+_{n,p}\} \Rightarrow C$. If $\{m,n\} \subset\nat,\ m, n \geq 2,\ m \not= n,$ then $RI^p_n(L) \not= RI^p_m(L).$ Further, if $p,q$ are  two distinct permutations, then $RI^p_n(L)
\not= RI^q_n(L).$ Hence, for each $n\in \nat,\ n \geq 2$, there exists $n!$ distinct general logic-systems that generate the same $C \in {\cal C}_f(L).$ Whether, for each $n\in \nat,\ n \geq 2,$ only one distinct permutation or each of the $n!$ permutations are utilized to define distinct general logic-systems, this implies that there exists a denumerable collection of general logic-systems each member of which generates $C$. \qed 

\noindent {\bf 4. GGU-model Operators.} \par\medskip
Of significance to physical science is the use of logic-systems to generate the development of a universe. For the General Grand Unification Model (GGU-model), logic-system behavior implies that physical-systems  are designed from rationally ordered combinations of constituents and each complete physical-system follows a rational development over observer-time. Their application to the GGU-model appears in Herrmann (2013a) and (2013b).\parm

\noindent{ \bf 5. A Formal Measurement of Intelligence.}\parm

General logic-systems can yield a measure for intelligence via the seventh Thurstone (1941) factor - ``Reasoning'' ability. For the GGU-model, the hyperfinite logic-system used is the $K^q_1$ as it is preserved by the operators ${IN}_q,\ G_q$ and St. Moreover, what follows is but one measure, among others, for the ability to reason. \parm

{\bf Definition 5.1} Intelligence, for GID-model, is the ability to apply rules specified by an algorithm and to obtain from a given logic-system distinct deductive conclusions or a specific conclusion. This ability is measured over a specific time interval. The measure itself is the number of reasoned distinct conclusions that can be obtained during that time interval or whether the final conclusion is the one specified.  \parm

Intelligence, as measured by Definition 5.1, has significant meaning via comparison. Consider the hyper-interval $\Hyper {[c_i, c_{i+1}]}$ and the hyperfinite logic-system $K^q_1(\lambda)$ restricted to this hyper-interval.  Consider the informal standard general logic-system $\r K^q_1$ obtained from $K^q_1$ by restriction.  Let agent $\r A$ be a standard agent that can perform only finitely many [i.e. $n$]  deductions over a time internal of length $c_{i+1} - c_i.$ (The first step is 
$\r F^q(\r t^q(i,0))$.) This is generalized to a set of ``superagents'' $\cal A$ where for each 
$n \in \nat,\ n>0,$ there is a member of $\cal A$ that can deduce $n$ distinct members of $\r d_q$ during this time interval.   Hence, for any $n \in \nat,\ n >0,$ there is a superagent $\r A_n$ that can obtain $n$ distinct deductions over time period $c_{i+1} - c_i.$ \pars

Formally characterizing the ``number'' of distinct deductions that a superagent can make, this number can be compared with hyperfinite set of deductions.  Consider the $\lambda$ in Theorems 4.q (Herrmann (2006b)). There exists a superagent agent $H$ that can deduce $\lambda +1$ distinct members of $d^q_x$. If one does not include the notion of superagents, then assume that an agent $H$ exists that can do hyper-deduction.  In mathematical logic, one can assign the superagent notion to such statements as ``for the formal predict logic and any 
$n \in \nat, n>0,$ there are well-formed formulas (formal theorems) that require $n$ or more steps to deduce.''  (There are multi-universe  models that do allow for superagents to exist in the sense that deductions can be continued via other agents indefinitely. Thus, in this case, a superagent is a finite collection of agents or, depending upon the cosmology,  a single agent.)   Definition 6.1 can be interpreted as follows: For an agent $H$ that can do hyper-deduction, agent $H$ is, in general, infinitely more intelligent than standard agent $\r A\in \cal A$ and, in general, can obtain conclusions that $\r A$ cannot. (In a few special cases, although it is not considered as deduction, special analysis can determine all the values of  $\{\Hyper {\b F}^q(\Hyper {\b t}^q(i,j))\mid 0\leq j \leq \lambda\}$.)\parm

\noindent{\bf 6. Potentially-Infinite.}\parm
    
This entire section has been removed since apparently the C-set theory axioms do not allow one to conclude that the set employed in Theorem 6.2 (i) exists.\parm

\centerline{\bf REFERENCES}\par\medskip

\id{B}ohm, D. 1957. ``Causality and Chance in Modern Physics,'' Harper \& Brothers, New York.\smallskip
\id{H}errmann, R. A. (2013a), ``Ultra-logic-systems applied to development paradigms,'' http://vixar.org/abs/1309.0004\hfil\break http://www.raherrmann.com/ultralsystem.pdf
\smallskip
\id{H}errmann, R. A. (2013b), ``Ultra-logic-systems applied to instruction paradigms,''  http://vixra.org/abs/1309.0125\hfil \break http://www.raherrmann.com/ultralsystem1.pdf
\smallskip
\id{H}errmann, R. A. (2006), ``General logic-systems and finite consequence operators,'' Logica Universalis 1:201-208. For a portion of this paper, see\hfil\break
http://arxiv.org/abs/math/0512559 \smallskip
\id{H}errmann, R. A. (2006a), {\it Logic for Everyone,}\hfil\break 
http://www.arxiv/abs/math/0601709\smallskip 
\id{H}errmann, R. A. (2006b), ``The GGU-model and generation of developmental paradigms,'' 
http://arxiv/org/abs/math/0605120\hfil\break
Latest verson http://vixra.org/abs/1308.0145\smallskip
\id{H}errmann, R. A. (2004), ``The best possible unification for any collection of physical theories,'' Internat. J. Math. and Math. Sci., 17:861-721.\hfil\break 
http://www.arxiv.org/abs/physics/0306147\smallskip
\id{H}errmann, R. A. (2004a), Nonstandard consequence operators generated by mixed logic-systems, \hfil\break http://arxiv.org/abs/math/0412562 
http://www.arxiv.org/abs/physics/0306147\smallskip
\id{H}errmann, R. A. (2002), ``Science Declares Our Universes IS Intelligently Designed,'' Xulon Press, Fairfax, VA.
\id{H}errmann, R. A. (2001), ``Hyperfinite and Standard Unifications for Physical Theories,'' Intern. J. Math. Math. Sci. {\bf 28}(2001), no.2, 93-102. http://arxiv.org/abs/physics/0105012 \smallskip
\id{H}errmann, Robert A. (2001a), ``Ultralogics and probability models,'' Internat. J. Math. and Math. Sci., 27(5):321-325. \hfil\break
http://www.arxiv/abs/physics/0105012\smallskip
\id{H}errmann, Robert A., (1999).``The encoding of quantum state information within subparticles.'' http://arxiv.org/abs/quant-ph/9909078 \smallskip
\id{H}errmann, Robert. A., (1979, 1993). The Theory of Ultralogics,  \hfil\break 
Part I at http://arxiv.org/abs/math/9903081 \hfil\break Part II at http://arxiv.org/abs/math/9903082 \smallskip
\id{H}errmann, R. A. (1987), ``Nonstandard Consequence operators,'' Kobe J. Math., 4(1):1-14. http://www.arxiv.org/abs/math.LO/9911204\smallskip
\id{J}ech, T. J., (1973), {\it The Axiom of Choice,} North-Holland, New York.\smallskip
 \id{J}ech, T. J., (1971), {\it Lectures in Set Theory with Particular Emphasis on the Method of Forcing,} No, 217, Lecture Notes is Mathematics, Springer-Verlag, New York. \smallskip 
\id{M}endelson, E., (1987), {\it Introduction to Mathematical Logic,} Wadsworth \& Brooks/Cole Advanced Book \& Software, Monterey, CA. \smallskip
\id{S}toll, Robert R. (1963), {\it Set Theory and Logic,} W. H. Freedom and Co, San Francisco, CA.\smallskip
\id{T}arski,  Alfred. (1956), {\it Logic, Semantics, Metamathematics; papers from 1923 - 1938},Oxford University Press, NY.
\smallskip
\id{W}ilder, R. L. (1967), {\it Introduction to The Foundations of Mathematics}, Wiley, NY.  
\smallskip
\id{W}\'ojcicki, R. (1981), ``On the content of logics part I. The representation theorem for lattices of logics,'' Reports on Mathematical Logic, 13:17-28.
\id{W}\'ojcicki, R. (1973), ``On matrix representations of consequence operators on \L ukasiewicz's Sentential Calculi,'' {\it Zeitschi. f. math. Logik und Grundlagen d. Math.,} 19:239-247.\medskip
\noindent {\it e-mail} drrangid@hotmail.com

\end